\documentclass[12pt,draftcls,onecolumn]{IEEEtran}

\usepackage{amsmath,amsfonts, amssymb,color}
\usepackage{epstopdf}
\usepackage{nomencl}
\usepackage{enumerate}
\usepackage{times} 
\usepackage{latexsym}
\usepackage{xspace}
\usepackage{amssymb,amsmath,amsfonts,color}
\usepackage{cite}
\usepackage{tikz,subcaption,cite}
\usetikzlibrary{decorations.pathmorphing}
\usetikzlibrary{shapes,arrows}
\usetikzlibrary{shadows.blur}
\usetikzlibrary{shapes.symbols}
\usepackage[american]{circuitikz}
\usetikzlibrary{circuits}

\tikzset{sin v source/.style={
  circle,
  draw,
  append after command={
    \pgfextra{
    \draw
      ($(\tikzlastnode.center)!0.5!(\tikzlastnode.west)$)
       arc[start angle=180,end angle=0,radius=0.425ex] 
      (\tikzlastnode.center)
       arc[start angle=180,end angle=360,radius=0.425ex]
      ($(\tikzlastnode.center)!0.5!(\tikzlastnode.east)$) 
    ;
    }
  },
  scale=1.5,
 }
}

\newenvironment{pfof}[1]{\vspace{1ex}\noindent{\itshape Proof of
    #1:}\hspace{0.5em}} {\hfill\oprocend\vspace{1ex}}

\newcommand\oprocendsymbol{\hbox{$\square$}}
\newcommand\oprocend{\relax\ifmmode\else\unskip\hfill\fi\oprocendsymbol}

\newtheorem{theorem}{Theorem}
\newtheorem{corollary}{Corollary}

\newtheorem{lemma}{Lemma}
\newtheorem{remark}{Remark}
\newtheorem{definition}{Definition}
\newtheorem{example}{Example}
\DeclareMathOperator{\trace}{trace}
\DeclareMathOperator{\diag}{diag}

\definecolor{cyan}{rgb}{0.0, 1.0, 1.0}

\begin{document}
\title{\LARGE \bf  System-Theoretic Performance Metrics for Low-Inertia Stability of Power Networks}

\author{Mohammad Pirani, John W. Simpson-Porco and  Baris Fidan 
\thanks{This material is based upon work supported in part by the Natural Sciences and Engineering Research Council of Canada (NSERC). M. Pirani and B. Fidan are with the Department of Mechanical and Mechatronics Engineering at the University of Waterloo, Waterloo, ON, Canada.  E-mail: \texttt{\{mpirani, fidan\}@uwaterloo.ca}. John W. Simpson-Porco is with the Department of Electrical and Computer Engineering at the University of Waterloo, Waterloo, ON, Canada. E-mail: {\texttt{jwsimpson@uwaterloo.ca}}. }
}

\maketitle

\thispagestyle{empty}
\pagestyle{empty}


\begin{abstract}

As bulk synchronous generators in the power grid are replaced by distributed generation interfaced through power electronics, inertia is removed from the system, prompting concerns over grid stability. Different metrics are available for quantifying grid stability and performance; however, no theoretical results are available comparing and contrasting these metrics. This paper presents a rigorous system-theoretic study of performance metrics for low-inertia stability.  For networks with uniform parameters, we derive explicit expressions for the eigenvalue damping ratios, and for the $\mathcal{H}_{2}$ and $\mathcal{H}_{\infty}$ norms of the linearized swing dynamics, from external power disturbances to different phase/frequency performance outputs.These expressions show the dependence of system performance on inertia constants, damping constants, and on the grid topology.
Surprisingly, we find that the $\mathcal{H}_2$ and $\mathcal{H}_{\infty}$ norms can display contradictory behavior as functions of the system inertia, indicating that low-inertia performance depends strongly on the chosen performance metric.

\end{abstract}


\section{Introduction}
\label{sec:intro}


Much attention has recently been focused on the integration of renewable energy sources into large-scale electric power systems. While traditional synchronous generators are characterized by large rotating inertias, renewables are typically integrated through power converters which are purely electronic and therefore provide no inertial response. As renewables supplant traditional generation, the total inertia present in the grid decreases, leading to concerns over ``low-inertia stability'' of such renewable-dominated systems \cite{Telegina}.

Quantifying the effects of lowered inertia on power grid stability, transients, and sensitivity to disturbances is a topic of present interest. In this direction, the effect of low rotational inertia on system stability was studied in \cite{Ulbig1, Ulbig2}, where effects were quantified in terms of (i) transients after a fault, and (ii) the region of attraction of a stable equilibrium point. The authors showed that grid topology can play a significant role when inertia levels are heterogeneous throughout the grid. In \cite{Borsche1}, the effects of lowered inertia on eigenvalue damping ratios and on frequency overshoot was studied, and an optimization problem was posed to determine optimal inertia values which maximize damping ratios while ensuring admissible transient behavior after a large disturbance.
%
%

Another method for quantifying power system performance is via a system norm, which measures the sensitivity of a chosen performance output to external disturbances. The $\mathcal{H}_2$ performance of the swing dynamics was studied in \cite{Tegling}, where phase differences of the network were the chosen performance outputs. Interestingly, in this case the norm was found to be independent of both network topology and inertia values. An optimal inertia placement problem for minimizing the system's $\mathcal{H}_2$ norm was introduced in \cite{Poola}, by considering frequency deviations and phase differences as output measurements. In \cite{Amer} the effect of disturbances on frequency deviations was studied by optimizing the $\mathcal{H}_2$ norm, the $\mathcal{H}_{\infty}$ norm, and the locations of eigenvalues. Sensitivity of the dominant eigenvalue to variations in inertia was considered in \cite{Ding}, and the zeros of swing dynamics with frequency outputs was studied in \cite{Koorehdavoudi}.

In summary, various metrics have been proposed for quantifying low-inertia stability. An important question to ask is whether these metrics are always consistent with one another. That is, if one metric shows a degradation in system performance, do the others? Unfortunately, we will show that the answer in general is no, and that these metrics can even yield contradictory results. 

Our approach is to analytically study the linearized swing dynamics of the network. We first consider the case of a single generator, the so-called single-machine infinite-bus (SMIB) system, and derive closed-form results for (i) the $\mathcal{H}_2$ and $\mathcal{H}_{\infty}$ norms, for phase  cohesiveness output, and (ii) the eigenvalues of the system. Surprisingly, for phase output we find that the $\mathcal{H}_{\infty}$ norm is an \emph{increasing} function of system inertia. In other words, the system becomes more robust as inertia is \emph{removed}. We then move to the case of a network of generators, and extend our single-machine results under the assumption of uniform inertia and damping coefficients \cite{Tegling}. In this network case, we show that the $\mathcal{H}_{\infty}$ norm depends on the algebraic connectivity $\lambda_2$ of the grid's admittance matrix. Our work can also be interpreted as a further contribution to the theory of robust networked dynamical systems \cite{Fitch2,Bamieh2,Jovanovic,Summers,Scardovi,Siami,ArxiveRobutness}.


%


The paper is organized as follows. Section \ref{sec:influence} describes the modeling of the power network and swing dynamics with phase cohesiveness and frequency performance outputs. 
In Section \ref{Sec:SMIB} we study the case of a single machine, pedagogically explaining our main results in terms of Bode plots and eigenvalues. Section \ref{sec:doubleeee} contains our main technical results, where we derive expressions for the $\mathcal{H}_2$ and $\mathcal{H}_{\infty}$ norms of the swing dynamics for each performance output, and discuss the dependence of each robustness metric on inertia, damping, and network connectivity. Finally we conclude in Section \ref{sec:conclusion}. The remainder of this section establishes some notation.

\subsection{Notation and Definitions}
\label{sec:not}

In this paper, an undirected network is denoted  by  $\mathcal{G}=\{\mathcal{V},\mathcal{E}\}$,  where $\mathcal{V} = \{1,\ldots,n\}$ is a set of nodes and $\mathcal{E} \subset \mathcal{V}\times\mathcal{V}$ is the set of edges.  Neighbors of node $i \in \mathcal{V}$ are given by the set $\mathcal{N}_i = \{j \in \mathcal{V} \mid (i, j) \in \mathcal{E}\}$. The adjacency matrix of the graph is the symmetric $n \times n$ matrix $A$, where $A_{ij} > 0$ if $(i,j) \in \mathcal{E}$ and zero otherwise. The degree of node $i$ is denoted by  $d_i \triangleq \sum_{j=1}^nA_{ij}$. The Laplacian matrix of the graph is given by $L \triangleq D - A$, where $D = \diag(d_1, d_2, \ldots, d_n)$.  The eigenvalues of the Laplacian are real and nonnegative, and are denoted by $0 = \lambda_1(L) \le \lambda_2(L) \le \ldots \le \lambda_n(L)$.  The $i$th eigenvalue of the Laplacian matrix is simply denote  by $\lambda_i$ in this paper, and we denote by $L^{\frac{1}{2}}$ the matrix square root of $L$.

\section{ Power Network Model}
\label{sec:influence}

Consider a power transmission network with $n$ buses $\mathcal{V}=\{1,\ldots,n\}$ and a set of transmission lines $\mathcal{E}$. Here we assume a Kron-reduced transmission network model, where all buses are modeled as generators and branch resistances are neglected \cite{Tegling}. At each bus $i \in \mathcal{V}$, there is a generator with inertia constant $M_i > 0$, damping/droop constant $D_i > 0$, and voltage phase angle $\theta_i$. The dynamics of the $i$th generator is described by the \emph{swing equation}
\begin{equation}
M_i\ddot{\theta}_i+D_i\dot{\theta}_i=P_{m,i}-P_{e,i}+w_i(t),
\label{eqn:mainpower}
\end{equation}
where $P_{m,i}$ is the constant mechanical power input from turbine and $w_i(t)$ models disturbances arising from  generation or  local load variations. The term $P_{e,i}$ is the real electrical power injected from $i$-th generator to the network, given by
\begin{align}
P_{e,i}= \sum_{j\in \mathcal{N}_i}\nolimits V_iV_jB_{ij}\sin(\theta_i-\theta_j),
\label{eqn:gednn}
\end{align}
where $V_i$ is the nodal voltage magnitude and $-B_{ij} < 0$ is the susceptance associated with edge $(i,j)\in\mathcal{E}$. We further approximate \eqref{eqn:gednn} using the so-called DC Power Flow, where $V_i \simeq V_j \simeq 1$ and $|\theta_i-\theta_j| <\!\!< 1$, leading to the linear model
\begin{align}
P_{e,i}\approx \sum_{j\in \mathcal{N}_i}\nolimits B_{ij} (\theta_i-\theta_j).
\label{eqn:linearized}
\end{align}
Substituting \eqref{eqn:linearized} into  \eqref{eqn:mainpower} yields
\begin{align}
M_i\ddot{\theta}_i+D_i\dot{\theta}_i\approx - \sum_{j\in \mathcal{N}_i}\nolimits B_{ij} (\theta_i-\theta_j)+P_{m,i}+w_i.
\label{eqn:genn}
\end{align}

In this paper we assume homogeneous inertia and damping parameters, i.e., $M_i=M$ and $D_i=D$ for all $i=1,\ldots,n$ similarly to \cite{EmmaTegling}. This assumption  allows us to establish closed-form expressions for our results. After shifting the equilibrium point of \eqref{eqn:genn} to the origin, the term proportional to $P_{m,i}$ may be removed and the dynamics of the generators can be written in state-space form for $\boldsymbol{\theta}=[\theta_1 , ..., \theta_n]^T$ and $\mathbf{w}=[w_1 , ..., w_n]^T$ as
\begin{align}
\begin{bmatrix}
       \dot{\boldsymbol{\theta}} \\[0.3em]
     \ddot{\boldsymbol{\theta}}
     \end{bmatrix}&=\underbrace{\begin{bmatrix}
       \mathbf{0}_{n} & I_{n}         \\[0.3em]
     -\frac{1}{M}L & -\frac{D}{M}I_{n}
     \end{bmatrix}}_{A}\underbrace{\begin{bmatrix}
       {\boldsymbol{\theta}} \\[0.3em]
     \dot{\boldsymbol{\theta}}
     \end{bmatrix}}_{\boldsymbol{\Theta}}
     +\underbrace{\begin{bmatrix}
       \mathbf{0}_{n}        \\[0.3em]
    \frac{1}{M}I_{n}   
     \end{bmatrix}}_{F}\mathbf{w}(t),\nonumber\\
     \mathbf{y}&=C\boldsymbol{\Theta}\,,
     \label{eqn:powernetwork}
\end{align}
where $L$ is the Laplacian matrix with weights $B_{ij}$, and the output matrix $C$ can take several forms. With the aim of measuring useful quantifies for assessing system performance, we consider the following outputs:
\begin{enumerate}
\item[(i)] \textbf{Phase Cohesiveness:} $\mathbf{y} =C\boldsymbol{\Theta}= L^{\frac{1}{2}}\boldsymbol{\theta}$. With this choice, 
$$
\mathbf{y}^{\sf T}\mathbf{y} = \boldsymbol{\theta}^{\sf T}L\boldsymbol{\theta} = \sum_{\{i,j\}\in\mathcal{E}}\nolimits B_{ij}(\theta_i-\theta_j)^2\,,
$$
which measures how tightly phase angles are clustered in the network. This output was proposed in \cite{Tegling} to measure resistive losses during transients using the $\mathcal{H}_2$ norm, and has been more broadly used in the network control literature \cite{Fdolfler, Fulin}. An alternative way of defining this performance output is to use any other output matrix $\tilde{C}$ in \eqref{eqn:powernetwork} such that $\tilde{C}^{\sf T}\tilde{C}=L$. In this case, as both $\mathcal{H}_2$ and $\mathcal{H}_{\infty}$ norms are functions of the spectrum of $G^*G=F^{\sf T}(s^*I-A)^{-\sf T}C^{\sf T}C(sI-A)^{-1}F$, identical results will be obtained as if one used $C = [L^{\frac{1}{2}} \hspace{1mm} \mathbf{0}]$. One such choice is
\begin{equation}
   \mathbf{y} =\diag({B}_{ij})^{\frac{1}{2}}\mathcal{B}^{\sf T}\boldsymbol{\theta},
   \label{eqn:alternateoutput}
\end{equation}
where $\mathcal{B} \in \mathbb{R}^{n \times |\mathcal{E}|}$ is the incidence matrix associated with the network. In this case we have an output associated with each edge. In fact, $y_{ij}=B_{ij}^{\frac{1}{2}}(\theta_i-\theta_j)$ which is proportional to the power transmitted across line $\{i,j\}$. Therefore, this output can be interpreted either as a measure of coherence (cohesiveness), power losses, or power flow on transmission lines.
\item[(ii)] \textbf{Frequency:} $\textbf{y} =C\boldsymbol{\Theta}= \dot{\boldsymbol{\theta}}$. Large frequency transients are unacceptable during operations, and therefore quantifying the effect of disturbances on frequency is important \cite{Amer}.
\item[(iii)] \textbf{Phase Cohesiveness \& Frequency:}  Combining the previous two outputs, we obtain
$$
\mathbf{y} =C\boldsymbol{\Theta}= \begin{bmatrix}
L^{\frac{1}{2}}\boldsymbol{\theta}\\
\kappa \dot{\boldsymbol{\theta}}
\end{bmatrix}\,,
$$
where $\kappa > 0$ is a design parameter. This performance output was used in \cite{Poola} in the context of optimizing the placement of inertia in the grid, and aims to simultaneously capture phase and frequency performance.
\end{enumerate}

The performance metrics we are investigating in this paper are (a) the poles of the swing dynamics \eqref{eqn:powernetwork} (eigenvalues of the $A$ matrix), which provide a stability measure independent of the chosen output, and (b) system $\mathcal{H}_2$ and $\mathcal{H}_{\infty}$ norms of \eqref{eqn:powernetwork}, defined as 
\begin{align}
||G||_2 &\triangleq \left( \frac{1}{2\pi}\trace\int_0^{\infty}G^*(j\omega)G(j\omega)d\omega \right)^{\frac{1}{2}},\nonumber \\
||G||_{\infty} &\triangleq \sup_{\omega\in \mathbb{R}}{\lambda_{max}^{\frac{1}{2}}(G^*(j\omega)G(j\omega))},
\end{align}
where $G(.)$ is the transfer function from external disturbance $\mathbf{w}(t)$ to different performance outputs mentioned above. In Section \ref{sec:doubleeee}, we derive closed-form expressions for the poles and damping ratios of  \eqref{eqn:powernetwork} and $\mathcal{H}_{2}$ and $\mathcal{H}_{\infty}$ norms for outputs (i) and (ii) mentioned before. These expressions are in terms of the spectrum of the Laplacian matrix as well as physical parameters of the system. Output (iii) proved too difficult to study analytically. However, we demonstrate numerically  that in general, the corresponding $\mathcal{H}_2$ and $\mathcal{H}_{\infty}$ norms for output (iii) show contradictory behaviour as a function of the inertia. 

\section{Low-Inertia Performance of Single-Machine Infinite-Bus (SMIB) System}
\label{Sec:SMIB}

\begin{figure}[h!]
\centering
\includegraphics[width=0.4\linewidth]{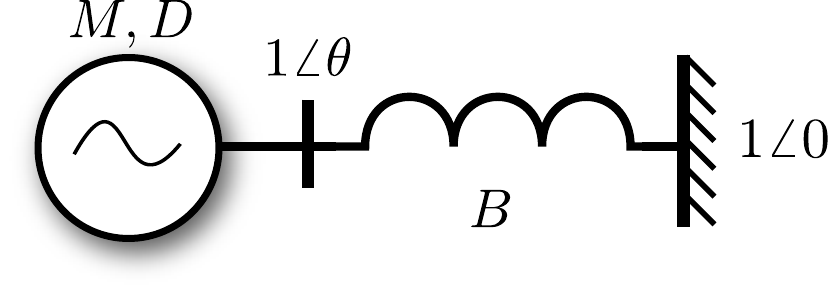}
\caption{A single generator with voltage phase angle $\theta$ connected to a large power system.}
\label{Fig:SMIB}
\end{figure}

Before proceeding to a more general setting consisting of many generators interacting over a network, we build intuition by considering the case of a single machine connected to a large power system (an ``infinite bus''), shown in Figure \ref{Fig:SMIB}.

The SMIB system is described by the linearized dynamics
\begin{equation}\label{Eq:SMIB}
\begin{aligned}
M\ddot{\theta} &= -D\dot{\theta} + P - B\theta + w\,,\\
y &= B^\frac{1}{2}\theta,
\end{aligned}
\end{equation}
where $M, D, B > 0$ and $P$ are real scalars. For this case, the output $y$ corresponds to the phase cohesiveness output described in Section \ref{sec:influence}. The following result follows as a special case of the more general result presented in Section \ref{sec:doubleeee}; frequency output results are deferred to the next section.

\begin{theorem}\label{Thm:SMIB}\textbf{($\mathcal{H}_2$ and $\mathcal{H}_{\infty}$ Performance of SMIB System)}:
Consider the single-machine infinite-bus system described by the dynamics \eqref{Eq:SMIB}, with the phase output $y = B^{\frac{1}{2}}\theta$. The $\mathcal{H}_2$ and $\mathcal{H}_{\infty}$ norms of the system are
 \begin{align}
||G||_{2}=\left(\frac{1}{2D}\right)^{\frac{1}{2}},
  \label{Eq:SMIB-H2}
\end{align}
and
 \begin{align}
||G||_{\infty}=\begin{cases}
  \frac{2M\sqrt{B}}{D\sqrt{4MB-D^2}},  & \quad \text{if } \frac{D^2}{2MB}\leq 1,\\
   \frac{1}{\sqrt{B}} & \quad  \text{otherwise}.\\
  \end{cases}
  \label{Eq:SMIB-HInf}
\end{align}
\end{theorem}

We are primarily interested in the parametric dependence of \eqref{Eq:SMIB-H2} and \eqref{Eq:SMIB-HInf} on the inertia constant $M$, and make two main observations.
First, the $\mathcal{H}_2$ norm \eqref{Eq:SMIB-H2} is independent of $M$. This indicates that the RMS or ``average'' sensitivity of the system to disturbances will be the same whether inertia is large or small. Second, the $\mathcal{H}_{\infty}$ norm \eqref{Eq:SMIB-HInf} is independent of $M$ for $M \in (0,D^2/2B)$, and strictly increasing in $M$ for $M \in [D^2/2B,\infty)$. In sharp contrast to conventional wisdom then, the system becomes more robust as inertia is removed.

\begin{figure}[h!]
\centering
\includegraphics[width=0.5\linewidth]{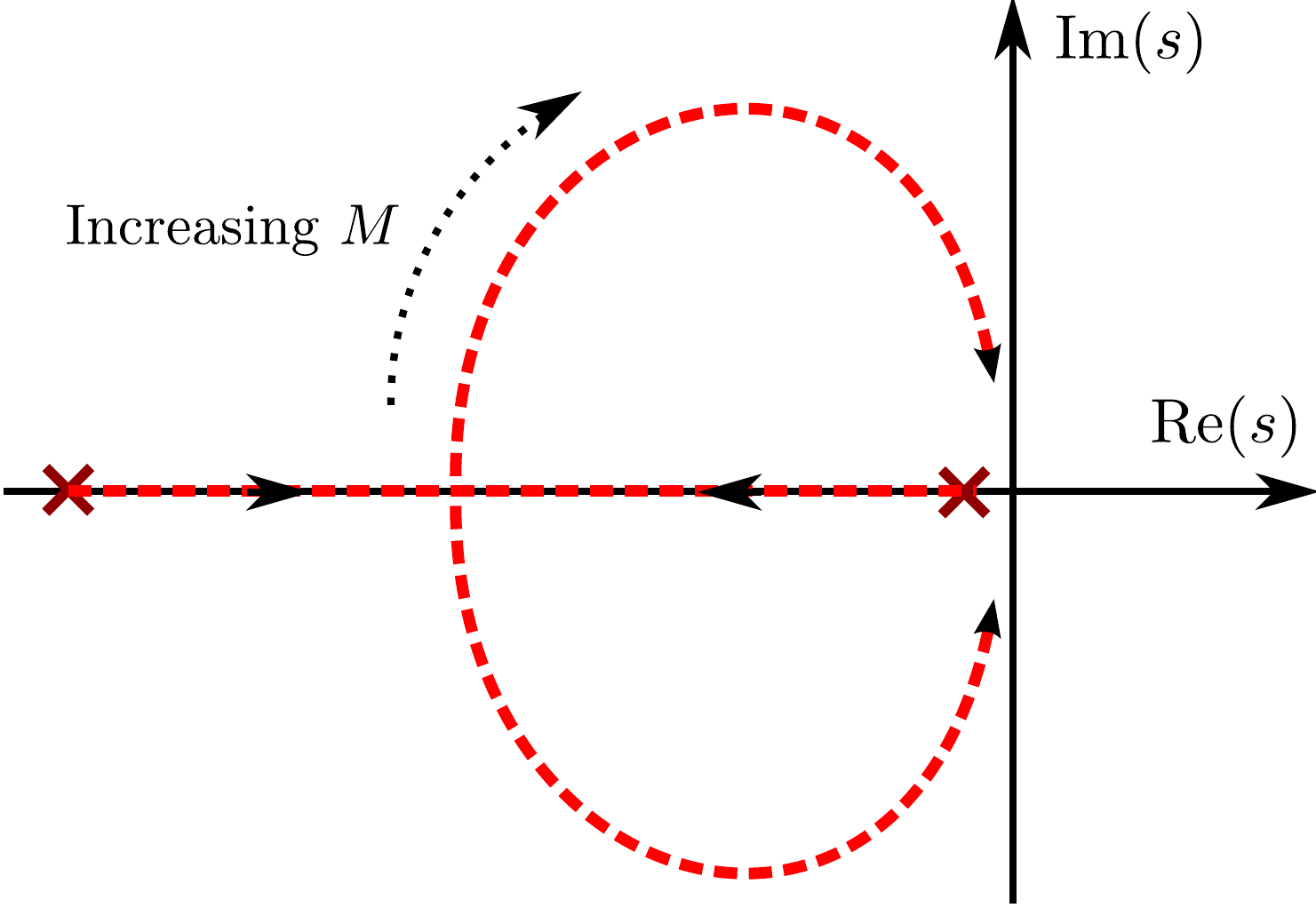}
\caption{ Root locus of SMIB system as a function of inertia.}
\label{Fig:SMIB-Rlocus}
\end{figure}

\begin{figure}[h!]
\centering
\includegraphics[width=0.7\linewidth]{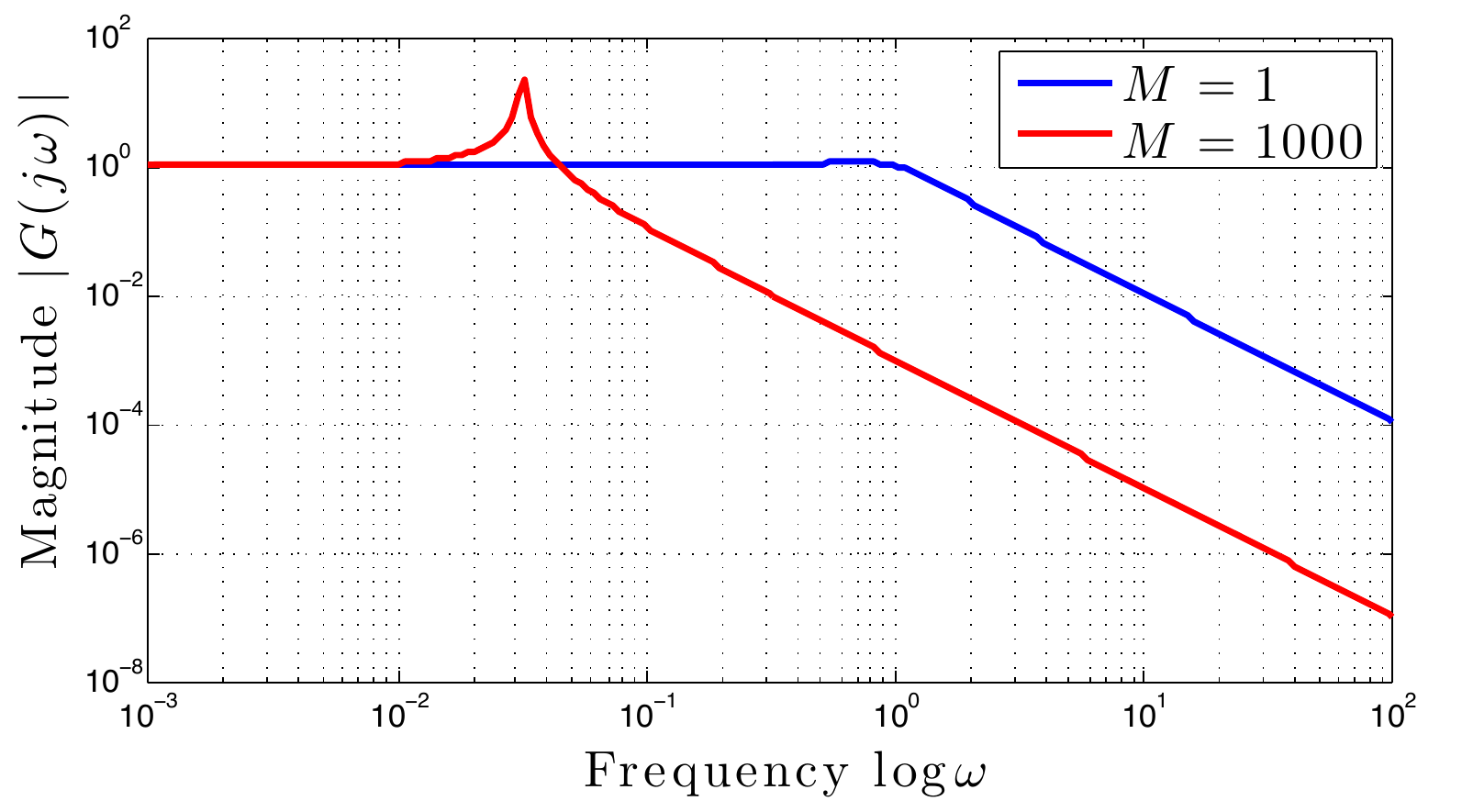}
\caption{Bode magnitude plot of SMIB system for two values of inertia.}
\label{Fig:SMIB-Bode}
\end{figure}

To understand this phenomena, consider the root locus plot (Figure \ref{Fig:SMIB-Rlocus}) of \eqref{Eq:SMIB} as a function of the inertia constant $M$.
The poles of the system \eqref{Eq:SMIB} are
$$
s = -\frac{D}{2M} \pm \frac{1}{2M}\sqrt{D^2-4MB}
$$
with natural frequency $\omega_n$ and damping ratio $\zeta$ given by
$$
\omega_n = \sqrt{B/M}\,,\qquad \zeta = \frac{D}{2\sqrt{BM}}\,.
$$
When $M$ is small, the system is heavily over-damped and the poles $s \in \{-D/M,-\epsilon\}$, where $0<\epsilon \ll D/M$ is a function of $M$, display a time-scale separation with fast and slow responses, respectively. As $M$ is increased, these poles converge on the real axis, break out into a complex conjugate pair, and eventually circle back to the origin. The damping ratio $\zeta$ continues to decrease however, as the poles converge faster to the imaginary axis than they do to the real axis. This results in an increasing peak in the Bode plot (Figure \ref{Fig:SMIB-Bode}), and therefore an increasing $\mathcal{H}_{\infty}$ norm. 

Conversely, the $\mathcal{H}_2$ result \eqref{Eq:SMIB-H2} indicates that despite this increasingly resonant peak in the Bode plot, the total (squared) area under the magnitude plot remains constant (Figure \ref{Fig:SMIB-Bode}). In summary, increasing the inertia makes the system increasingly resonant at the resonant frequency
$$
\omega_{\rm peak} = \omega_n\sqrt{1-2\zeta^2} = \sqrt{\frac{B}{M}\left(1-\frac{D^2}{2BM}\right)} \simeq \sqrt{\frac{B}{M}}\,,	
$$
while the magnitude roll-off occurs shortly after this resonant peak due to increased low-pass filtering from the large inertia. These results indicate that the relationship between inertia and system performance can be subtle, and depends strongly on the way performance is measured.

\section{$\mathcal{H}_{2}$ and $\mathcal{H}_{\infty}$ Robustness of the Swing Equation}
\label{sec:doubleeee}

This section contains our main technical results, extending the arguments from the SMIB system to a class of networks with homogeneous inertia and damping constants.


\subsection{Eigenvalues of the swing dynamics}

Our first result characterizes the eigenvalues (poles) of the linearized swing dynamics \eqref{eqn:powernetwork}.

\smallskip


\begin{theorem}\label{thm:dampratio}
Consider the power network described by the linearized swing dynamics \eqref{eqn:powernetwork}. The eigenvalues of \eqref{eqn:powernetwork} are given by
\begin{equation}\label{eqn:poledamp}
\begin{aligned}
s_{i1}&=-\frac{D}{2M}+ \frac{1}{2M}\sqrt{D^2-4M\lambda_i}\\
s_{i2}&=-\frac{D}{2M}- \frac{1}{2M}\sqrt{D^2-4M\lambda_i}
\end{aligned}\quad i=1,\ldots,n\,,
\end{equation}
and the smallest damping ratio $\zeta_{\rm min}$ of any mode equals
$$
\zeta_{\rm min} = \frac{D}{2\sqrt{M\lambda_n}}\,.
$$
\end{theorem}

\begin{IEEEproof}
The eigenvalues of $A$ are determined by $\mathrm{det}(sI_n-A) = 0$, which yields
\begin{align*}
{\rm det}\begin{bmatrix}
       sI_n & -I_{n} \\[0.3em]
     \frac{1}{M}L & sI_n+\frac{D}{M}I_{n}
     \end{bmatrix}
     &={\rm det} \left((s^2+s\frac{D}{M})I_n+\frac{1}{M}L\right)\\
     &=\prod_{i=1}^n\left(s^2+s\frac{D}{M}+\frac{1}{M}\lambda_i\right)=0\,,
\end{align*}
from which the expressions \eqref{eqn:poledamp} follow. By solving the pair of equations $2\zeta_i \omega_{n,i}=\frac{D}{M}$ and $\omega_{n,i}^2=\frac{\lambda_i}{M}$, the damping ratio of the $i$th mode is $\zeta_{i} = D/(2\sqrt{M\lambda_i})$ which obtains the result.
\end{IEEEproof}

\smallskip

While increasing the damping constant $D$ obviously damps the dynamics, Theorem \ref{thm:dampratio} indicates that, counter-intuitively, increasing inertia $M$ yields a less damped response. Moreover, the result shows that the largest eigenvalue $\lambda_n$ of the Laplacian matrix $L$ controls this minimally-damped mode.

\subsection{System norms for phase cohesiveness output}

We now present closed-form expressions for $\mathcal{H}_{2}$ and  $\mathcal{H}_{\infty}$  system norms of the swing dynamics \eqref{eqn:powernetwork}, from external disturbances $\mathbf{w}(t)$ to the phase cohesiveness performance output. The proof of case (i) of Theorem \ref{thm:sqrtc2} is presented in  \cite{Tegling} and the proof of case (ii) is in Appendix \ref{sec:appa}.

\begin{theorem}\label{thm:sqrtc2}\textbf{(Performance of Swing Dynamics with Phase Output)}:
Consider the power network described by the linearized swing dynamics \eqref{eqn:powernetwork} with the phase cohesiveness output $\mathbf{y}={L}^{\frac{1}{2}}\boldsymbol{\theta}$ or $\mathbf{y} = \diag({B}_{ij})^{\frac{1}{2}}\mathcal{B}^{\sf T}\boldsymbol{\theta}$. 
\begin{enumerate}
 
 \item[(i)] The $\mathcal{H}_{2}$ norm from disturbances to the output is
 \begin{align}
||G||_{2}=\left(\frac{n}{2D}\right)^{\frac{1}{2}},
  \label{eqn:cessj}
\end{align}

\item[(ii)] The $\mathcal{H}_{\infty}$ norm from disturbances to the output is
 \begin{align}
||G||_{\infty}=\begin{cases}
  \frac{2M\sqrt{\lambda_2}}{D\sqrt{4M\lambda_2-D^2}},  & \quad \text{if } \frac{D^2}{2M\lambda_2}\leq 1,\\
   \frac{1}{\sqrt{\lambda_2}} & \quad  \text{otherwise}.\\
  \end{cases}
  \label{eqn:cessss}
\end{align}
\end{enumerate}
\end{theorem}

In Corollary \ref{cor:behaviour}, we discuss the dependencies of $\mathcal{H}_{2}$ and $\mathcal{H}_{\infty}$ norms to system parameters. The proof is similar to that of Proposition 3 in \cite{PiraniSimpsonIFAC} and is omitted due to space limitations. 

\begin{corollary}
 System $\mathcal{H}_{2}$ norm  for phase cohesiveness output \eqref{eqn:cessj} is independent of the inertia $M$ and it is a monotonic decreasing function of the damping constant $D$. Moreover, the system $\mathcal{H}_{\infty}$ norm \eqref{eqn:cessss} is a continuously differentiable  and non-decreasing function of the inertia  $M$, and it is bounded from below as 
    \begin{equation}
    ||G||_{\infty}\geq \frac{1}{\sqrt{\lambda_2}},
    \label{eqn:bounded2}
    \end{equation}
    with strict equality sign for all $M\leq \frac{D^2}{2\lambda_2}$. Moreover, \eqref{eqn:cessss} is a convex function of $M$ for $M\leq \frac{D^2}{\lambda_2}$ and concave for $M> \frac{D^2}{\lambda_2}$. Furthermore, the $\mathcal{H}_{\infty}$ norm is a non-increasing function of $D$ and bounded from below by \eqref{eqn:bounded2}.
    \label{cor:behaviour}
\end{corollary}

\begin{figure}[h!]
\centering
\includegraphics[width=0.9\linewidth]{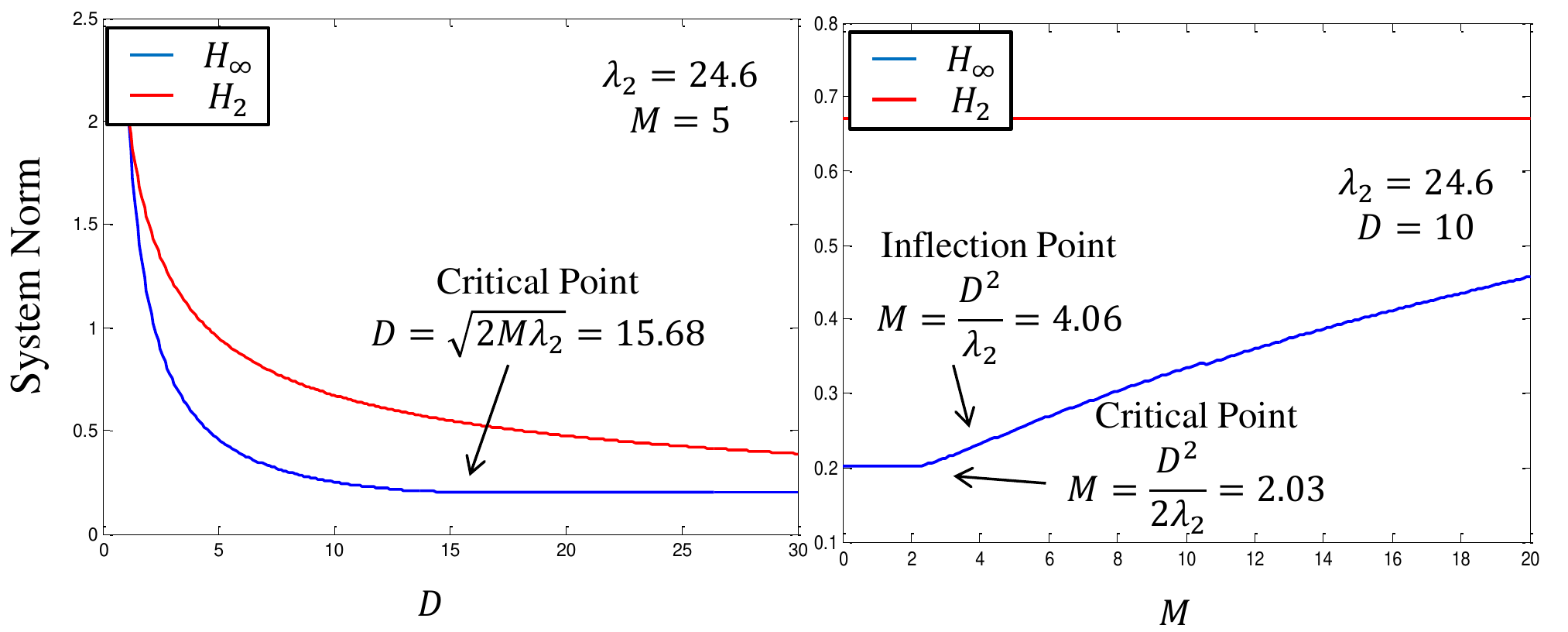}
\caption{System norms as functions of inertia and damping, for phase cohesiveness output.}
\label{fig:hinfpoasdeswef} 
\end{figure}

Fig.~\ref{fig:hinfpoasdeswef} shows the behavior of $\mathcal{H}_{2}$ and $\mathcal{H}_{\infty}$ norms of the swing dynamics \eqref{eqn:powernetwork} for phase cohesiveness output, as  functions of inertia $M$ and damping $D$.  As it is shown in Fig.~\ref{fig:hinfpoasdeswef} (left) and predicted by Corollary \ref{cor:behaviour}, system $\mathcal{H}_2$ norm is a monotonic decreasing function of $D$ and system $\mathcal{H}_{\infty}$ norm is monotonic decreasing function for $D\leq  \sqrt{2M\lambda_2}$ and is independent of $D$ for $D> \sqrt{2M\lambda_2}$. From Fig.~\ref{fig:hinfpoasdeswef} (right), the system $\mathcal{H}_{2}$ norm is independent of the inertia $M$ while the $\mathcal{H}_{\infty}$ norm is independent of $M$ for $M < \frac{D^2}{2\lambda_2}$ and increases by $M$ when $M \geq  \frac{D^2}{2\lambda_2}$ and changes its convexity at $M=\frac{D^2}{\lambda_2}$. 

\subsection{System norms for frequency output}

We now present closed-form expressions for $\mathcal{H}_{2}$ and $\mathcal{H}_{\infty}$  system norms of the swing dynamics \eqref{eqn:powernetwork}, from external disturbances $\mathbf{w}(t)$ to the frequency output. 

\smallskip

\begin{theorem}\label{thm:sqrdat2}\textbf{(Performance of Swing Dynamics with Frequency Output)}: Consider the power network described by the linearized swing dynamics \eqref{eqn:powernetwork} with frequency deviation output $\mathbf{y}=\dot{\boldsymbol{\theta}}$. 
\begin{enumerate}

\item[(i)] The $\mathcal{H}_{2}$ norm from disturbances to the output is
 \begin{align}
||G||_{2}=\left(\frac{n}{2DM}\right)^{\frac{1}{2}}.
\label{eqn:bsdfhb}
\end{align}

 \item[(ii)] The $\mathcal{H}_{\infty}$ norm from disturbances to the output is
 \begin{align}
||G||_{\infty}=\frac{1}{D}.
\label{eqn:frqout}
\end{align}
\end{enumerate}

\end{theorem}
\begin{IEEEproof}
 For case (i), we compute the $\mathcal{H}_2$ using the trace formula $||G||_2^2=\trace(F^{\sf T}PF)$, where $P$ is the observability Gramian $P=\int_{0}^{\infty}e^{A^{\sf T}t}C^{\sf T}Ce^{At}$ and it is uniquely obtained from the Lyapunov equation $P{A}+{A}^{\sf T}P=-C^{\sf T}C$. Here matrix $A$ is marginally stable and $(A,C)$ is not observable. However, since the mode corresponding to the marginally stable eigenvalue, $v=[\mathbf{1}_n^{\sf T} \hspace{1mm} \mathbf{0}_n^{\sf T}]^{\sf T}$ is not observable, i.e., $Ce^{At}v=Cv=\mathbf{0}_{2n}$ for all $t \geq 0$, and the rest of the eigenvalues are stable, the indefinite integral exists \cite{Doyle}. To calculate the observability Gramian, we have 
\begin{align}
\begin{bmatrix}
      P_{11} & P_{12}        \\[0.3em]
     P_{21} & P_{22}
     \end{bmatrix}A
     &+A^{\sf T}\begin{bmatrix}
      P_{11} & P_{12}        \\[0.3em]
     P_{21} & P_{22}
     \end{bmatrix} =\begin{bmatrix}
     \mathbf{0}_{n}  & \mathbf{0}_{n}       \\[0.3em]
     \mathbf{0}_{n} & -I_{n} 
     \end{bmatrix},
     \label{eqn:lyapp2}
\end{align}
Since $F= [\mathbf{0}_{n}, \frac{1}{M}I_{n}]^{\sf T}$, we have $F^{\sf T}PF=\frac{1}{M^2}P_{22}$; thus we only need to calculate $P_{22}$.  By solving \eqref{eqn:lyapp2} for $P_{22}$ we get $P_{22}=\frac{M}{2D}I_n$. Hence we have $||G||_2^2=\trace(F^{\sf T}PF)=\frac{n}{2DM}$. 
The proof of case (ii) is similar to case (ii) of Theorem \ref{thm:sqrtc2}. 
\end{IEEEproof}

The following corollary discusses the dependencies of system $\mathcal{H}_{2}$ and $\mathcal{H}_{\infty}$ norms \eqref{eqn:bsdfhb} and \eqref{eqn:frqout} to system parameters, inertia and damping constants. 

\begin{corollary}
 System $\mathcal{H}_{2}$ norm  for frequency output \eqref{eqn:bsdfhb} is a monotonic decreasing function of the inertia $M$ and the damping $D$. The $\mathcal{H}_{\infty}$ norm of the power network \eqref{eqn:frqout} is  an independent function of inertia and it is a monotonic decreasing function of $D$.
 \label{cor:dependee}
\end{corollary}

\begin{figure}[h!]
\centering
\includegraphics[width=0.9\linewidth]{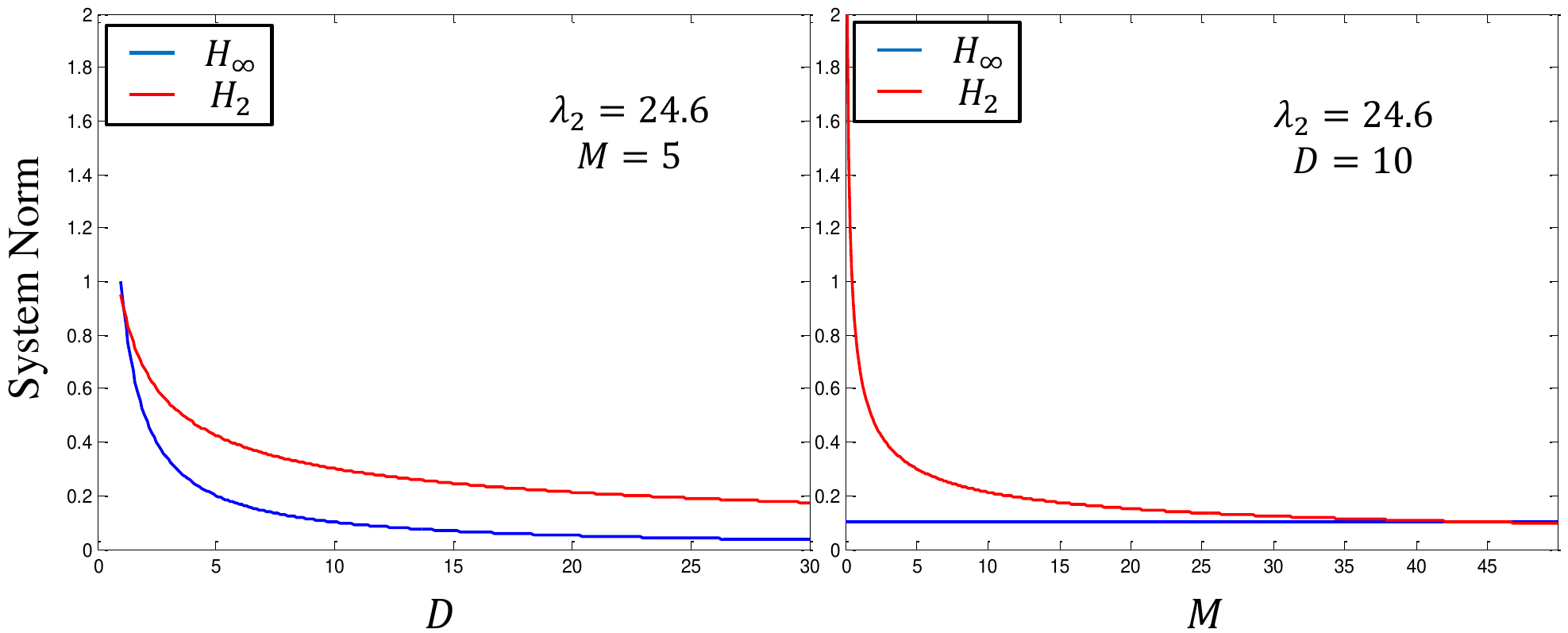}
\caption{System norms as functions of inertia and damping, for frequency output.}
\label{fig:hinfpoeswef} 
\end{figure}

Fig.~\ref{fig:hinfpoeswef} shows the behavior of $\mathcal{H}_{2}$ and $\mathcal{H}_{\infty}$ norms of the linearized swing dynamics \eqref{eqn:powernetwork} for frequency output, as  functions of inertia and damping. As it is shown in Fig.~\ref{fig:hinfpoeswef} and predicted by Corollary \ref{cor:dependee}, both metrics are monotonic decreasing functions of damping $M$ and inertia $D$ and the only exception is the invariance  of $\mathcal{H}_{\infty}$ with respect to variations of $M$, confirming \eqref{eqn:frqout}.

\begin{remark}\textbf{(Dependence on the Network Structure)}:
As it can be concluded from Theorem \ref{thm:sqrdat2}, both system $\mathcal{H}_2$ and $\mathcal{H}_{\infty}$ norms of \eqref{eqn:powernetwork} for frequency output case, are independent of the network structure. Such independence of network structure also holds for the system $\mathcal{H}_2$ norm for phase  output \eqref{eqn:cessj}, based on Theorem  \ref{thm:sqrtc2}. However, for this particular performance output, system $\mathcal{H}_{\infty}$ norm \eqref{eqn:cessss} is highly dependent on the connectivity of the underlying network. 
\end{remark}

\subsection{Combined phase cohesiveness and frequency outputs}
\label{Sec:PhaseFrequency}

\begin{figure}[h!]
\centering
\includegraphics[width=0.6\linewidth]{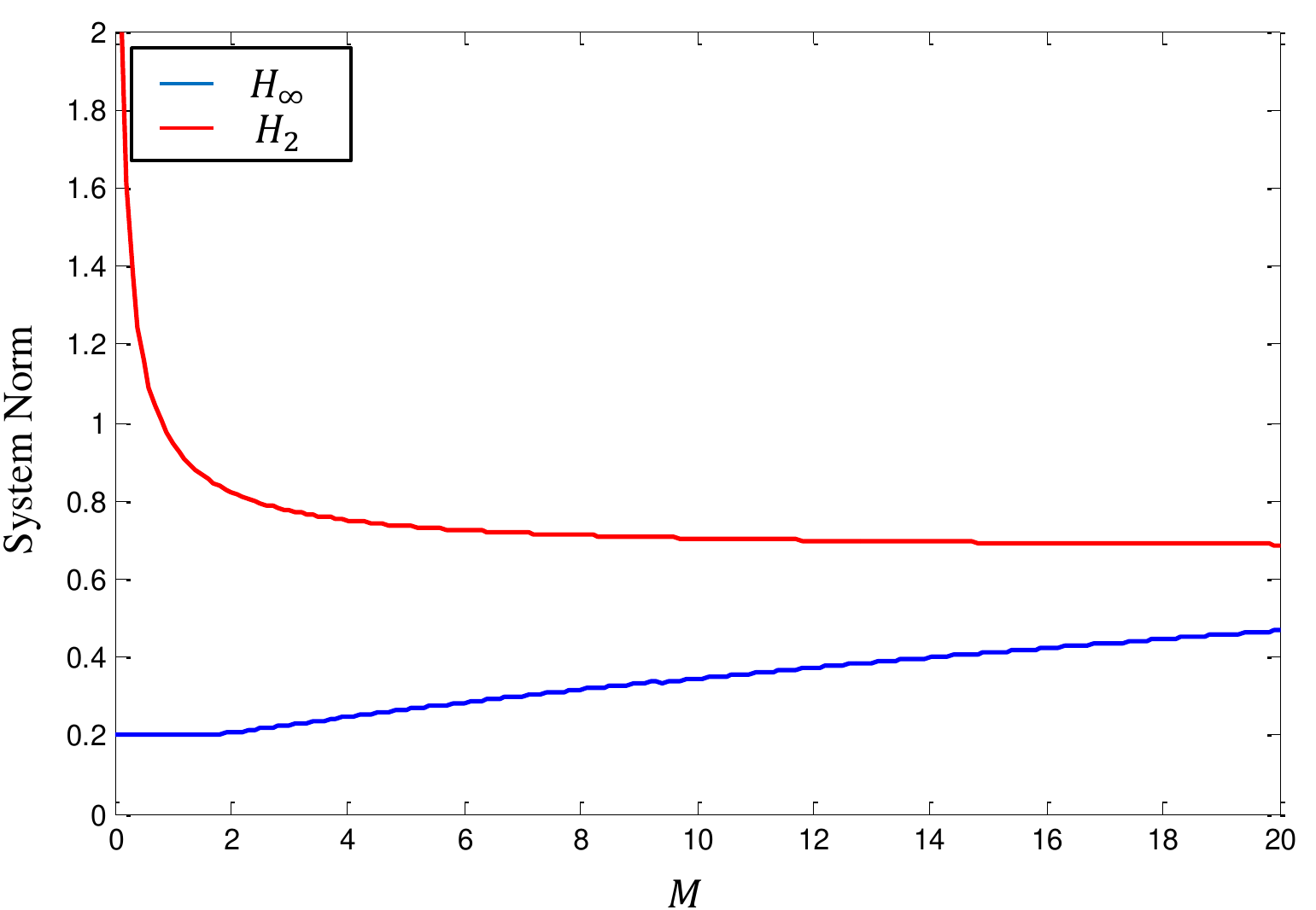}
\caption{System norms as functions of inertia for output \eqref{eqn:poolaa}.}
\label{fig:Poola}
\end{figure}

Finally, we consider the output proposed in \cite{Poola} which simultaneously accounts for phase cohesiveness and frequency deviations:
\begin{equation}
    \mathbf{y}=\begin{bmatrix}
       L^{\frac{1}{2}} & \mathbf{0}_{n}         \\[0.3em]
     \mathbf{0}_{n} & \kappa I_{n}
     \end{bmatrix} \begin{bmatrix}
       \boldsymbol{\theta}      \\[0.3em]
     \dot{\boldsymbol{\theta}} 
     \end{bmatrix},
     \label{eqn:poolaa}
\end{equation}
where $\kappa > 0$ is a chosen constant. Intuitively, based on results from Theorem \ref{thm:sqrtc2} and Theorem \ref{thm:sqrdat2} we expect that with the output \eqref{eqn:poolaa} (i) the $\mathcal{H}_{\infty}$ should be an increasing function of inertia, and (ii) the $\mathcal{H}_2$ norm should be a decreasing function of inertia. Figure \ref{fig:Poola} shows the trace of both system norms obtained numerically.

\section{ Conclusions}
\label{sec:conclusion}

In this paper we studied various metrics for quantifying performance in low-inertia power systems. Under the assumption of homogeneous inertia and damping parameters, we derived a closed-form expression for the minimally damped eigenvalue, and calculated the $\mathcal{H}_2$ and $\mathcal{H}_{\infty}$ system norms for phase cohesiveness and frequency deviation outputs. These expressions depend on the network structure through the spectrum of the Laplacian/admittance matrix. Our results show that these various metrics of performance do not necessarily trend in the same direction as a function of grid inertia; in general, they are competing objectives.
As the derived system norms are functions of both physical and network properties, optimizing these system norms with respect to either the physical or network structure is an important field of future research. Another avenue for extending the results presented in this paper is to quantify these system norms when system parameters are no longer homogeneous.


\bibliographystyle{IEEEtran}
\bibliography{main}

\begin{appendix}

\subsection{Proof of Theorem \ref{thm:sqrtc2}}\label{sec:appa}
We use the following lemma to prove Theorem \ref{thm:sqrtc2}, case (ii).
\begin{lemma}\label{Lem:HinfTransform}
Consider the square linear system
\begin{equation}\label{Eq:ControlSystem}
G:\,\begin{cases}
\begin{aligned}
\dot{x} &= Ax + Fu\\
y &= Cx\,.
\end{aligned}
\end{cases}
\end{equation}
with $x \in \mathbb{R}^n$ and $u,y \in \mathbb{R}^m$, and for orthogonal $V \in \mathbb{R}^{m \times m}$ consider the input/output transformation $\tilde{y} = Vy$, $\tilde{u} = Vu$, leading to the linear system
\begin{equation}\label{Eq:ControlSystem2}
\widetilde{G}:\,\begin{cases}
\begin{aligned}
\dot{x} &= Ax + FV^{-1}\tilde{u}\\
\tilde{y} &= VCx\,.
\end{aligned}
\end{cases}
\end{equation}
Then $\|G\|_{\infty} = \|\widetilde{G}\|_{\infty}$.
\label{lem:lemmaappendix}
\end{lemma}

\begin{IEEEproof}
The corresponding transfer functions are
\begin{align*}
G(s) &= C(sI-A)^{-1}F\\
\widetilde{G}(s) &= VC(sI-A)^{-1}FV^{-1}\,,
\end{align*}
and therefore
\begin{align*}
\widetilde{G}\widetilde{G}^{\star} &= VC(sI-A)^{-1}F\underbrace{V^{-1}V^{-{\sf T}}}_{=(V^{\sf T}V)^{-1} = I}F^{\sf T}(s^{\star}I-A)^{-{\sf T}}C^{\sf T}V^{\sf T}\\
&= V(GG^{\star})V^{\sf T} = V(GG^{\star})V^{-1}\,\,,
\end{align*}
where we have twice used that $V$ is orthogonal. Therefore, $GG^{\star}$ and $\widetilde{G}\widetilde{G}^{\star}$ are similar. It follows that for all $\omega \geq 0$
\begin{align*}
\sigma_{\rm max}(G(j\omega)) &= \lambda_{\rm max}(G(j\omega)G^{\sf T}(-j\omega))\\
&= \lambda_{\rm max}(\widetilde{G}(j\omega)\widetilde{G}^{\sf T}(-j\omega))\\
&= \sigma_{\rm max}(\widetilde{G}(j\omega))\,,
\end{align*}
and the result follows by taking supremums over $\omega$.
\end{IEEEproof}

\begin{pfof}{Theorem \ref{thm:sqrtc2}}
The model \eqref{eqn:powernetwork} has state-space matrices
$$
A = \begin{bmatrix}
0 & I_n\\ -\frac{1}{M}L & -\frac{D}{M}I
\end{bmatrix}\,,\quad F = \begin{bmatrix}
0 \\ \frac{1}{M}I_n
\end{bmatrix}\,, \quad C = \begin{bmatrix}
L^{\frac{1}{2}} & 0
\end{bmatrix}\,,
$$
with state vector $\boldsymbol{\Theta} = (\boldsymbol{\theta},\dot{\boldsymbol{\theta}})$. Let $\Lambda = V^{\sf T}LV$ be the eigendecomposition of $L$, where $V$ may be taken to be orthogonal. Consider the invertible change of states $\tilde{\boldsymbol{\Theta}} = (V^{\sf T}\boldsymbol{\theta},V^{\sf T}\dot{\boldsymbol{\theta}})$. Then a straightforward computation shows that
\begin{equation}\label{Eq:SwingTransformed1}
\begin{aligned}
\dot{\tilde{\boldsymbol{\Theta}}} &= \begin{bmatrix}
0 & I_n\\ -\frac{1}{M}\Lambda & -\frac{D}{M}I
\end{bmatrix}\tilde{\boldsymbol{\Theta}} + \begin{bmatrix}
0 \\ \frac{1}{M}V^{\sf T}
\end{bmatrix}w\\
y &= \begin{bmatrix}
L^{\frac{1}{2}}V & 0
\end{bmatrix}\tilde{\boldsymbol{\Theta}}\,.
\end{aligned}
\end{equation}
The model \eqref{Eq:SwingTransformed1} has the same transfer function as \eqref{eqn:powernetwork}, and hence the same system norm. Now consider an input/output transformation on \eqref{Eq:SwingTransformed1}, where $\bar{y} = V^{\sf T}y$ and $\bar{w} = V^{\sf T}w$\,. Then by Lemma \ref{Lem:HinfTransform}, the transformed system
\begin{equation}\label{Eq:SwingTransformed2}
\begin{aligned}
\dot{\tilde{\boldsymbol{\Theta}}} &= \begin{bmatrix}
0 & I_n\\ -\frac{1}{M}\Lambda & -\frac{D}{M}I
\end{bmatrix}\tilde{\boldsymbol{\Theta}} + \begin{bmatrix}
0 \\ \frac{1}{M}\underbrace{V^{\sf T}V}_{= I_n}
\end{bmatrix}\bar{w}\\
\bar{y} &= 
\underbrace{
\begin{bmatrix}
V^{\sf T}L^{\frac{1}{2}}V & 0
\end{bmatrix}}_{=\begin{bmatrix}\Lambda^{\frac{1}{2}} & 0 \end{bmatrix}}
\tilde{\boldsymbol{\Theta}}\,.
\end{aligned}
\end{equation}
has the same system norm as \eqref{Eq:SwingTransformed1}. The system \eqref{Eq:SwingTransformed2} is comprised of $n$ decoupled subsystems, each of the form
\begin{equation}\label{Eq:SwingTransformed3}
\begin{aligned}
\dot{\tilde{\boldsymbol{\Theta}}}_i &= \begin{bmatrix}
0 & 1\\ -\frac{1}{M}\lambda_i & -\frac{D}{M}
\end{bmatrix}\tilde{\boldsymbol{\Theta}}_i + \begin{bmatrix}
0 \\ \frac{1}{M}
\end{bmatrix}\bar{w}_i\\
\bar{y}_i &= 
\begin{bmatrix}
\lambda_i^{\frac{1}{2}} & 0
\end{bmatrix}
\tilde{\boldsymbol{\Theta}}_i\,.
\end{aligned}
\end{equation}
with transfer functions
$$
\tilde{G}_i(s) = \frac{\lambda_i^{\frac{1}{2}}}{Ms^2 + Ds + \lambda_i}\,, \qquad i \in \{1,\ldots,n\}\,.
$$
Clearly $\tilde{G}_1(s) = 0$. For $i \in \{2,\ldots,n\}$,  we have 
\begin{align*}
|\tilde{G}_i(j\omega)|^2&=\tilde{G}_i(-j\omega)\tilde{G}_i(j\omega)=\frac{\lambda_i}{\underbrace{(\lambda_i-M\omega^2)^2+D^2\omega^2}_{f(\omega)}}.
\end{align*}
Maximizing $|\tilde{G}_i(j\omega)|^2$ with respect to $\omega$ is equivalent to minimizing $f(\omega)$. By setting $\frac{df(\omega)}{d\omega}=0$ we get $\bar{\omega}_1=0$ and $\bar{\omega}_2=(\frac{\lambda_i}{M}-\frac{D^2}{2M^2})^{\frac{1}{2}}$ as critical points. Here 
$\bar{\omega}_2$ is the global minimizer of $f(\omega)$, unless $\frac{D^2}{2M\lambda_i}>1$. Substituting these critical values back into the formula for $|\tilde{G}_i(j\omega)|^2$, we find for $i \in \{2,\ldots,n\}$ that
\begin{align}
||\tilde{G}_i||_{\infty}=\begin{cases}
  \frac{2M\sqrt{\lambda_i}}{D\sqrt{4M\lambda_i-D^2}},  & \quad \text{if } \frac{D^2}{2M\lambda_i}\leq 1,\\
   \frac{1}{\sqrt{\lambda_i}} & \quad  \text{otherwise}\,.\\
  \end{cases}
  \label{Eq:HinfofIndividualSubsystems}
\end{align}
Since $0 < \lambda_2 < \lambda_2 \leq \lambda_3 \leq \cdots \leq \lambda_n$ and $||\tilde{G}_i||_{\infty}$ is a monotonically decreasing function of $\lambda_i$,  the result follows. 
\end{pfof}

\end{appendix}

\end{document}